\documentclass[leqno,b5paper]{amsart}

\usepackage{tikz}
\usepackage{tikz-cd}
\usepackage{blindtext}
\usetikzlibrary{automata, arrows.meta, positioning, decorations.markings}
\usepackage{amsmath,amsfonts,amsthm,latexsym}
\usepackage{mathrsfs, textcomp}
\newtheorem{theorem}{Theorem}[section]
\newtheorem{lemma}[theorem]{Lemma}
\newtheorem{proposition}[theorem]{Proposition}
\newtheorem{corollary}[theorem]{Corollary}

\theoremstyle{definition}
\newtheorem{definition}[theorem]{Definition}
\newtheorem{example}[theorem]{Example}

\numberwithin{equation}{section}

\usepackage{amssymb}
\usepackage{amscd}
\usepackage{eqlist, color}
\usepackage{mathrsfs, textcomp}
\usepackage{color, eqlist}
\usepackage{graphicx}
\usepackage{bm}
\usepackage{amsmath,amsfonts,amsthm,latexsym}
\usepackage[all]{xy}
\usepackage{tikz}
\usepackage{tikz-cd}

\DeclareMathOperator{\inte}{int}

\DeclareMathOperator{\Int}{int}

\DeclareMathOperator{\cl}{cl}

\theoremstyle{definition}\newtheorem{thm}{Theorem}[section]
\theoremstyle{definition}
\theoremstyle{definition}
\theoremstyle{definition}
\theoremstyle{definition}
\theoremstyle{remark}
\theoremstyle{definition}

\theoremstyle{definition}
\newtheorem{observation}[thm]{Observation}

\begin{document}
	\title[On the category of $(i,j)$-Baire bilocales]
	{On the category of $(i,j)$-Baire bilocales}

	\author{Mbekezeli Nxumalo}
	\address{Department of Mathematics, Rhodes University, P.O. Box 94, Grahamstown 6140, South Africa.}
	\email{sibahlezwide@gmail.com}
	%\email{bloses@unisa.ac.za}
	\subjclass[2010]{06D22}
	\keywords {$(i,j)$-Baire, topobilocale, $i$-prefit, $i$-pseudocomplete, ideal bilocale, relatively $(i,j)$-Baire}
	%\thanks{This paper was communicated by Prof. Themba Dube.}
	
	%%% ----------------------------------------------------------------------
	
	\let\thefootnote\relax\footnote{}
	
	\begin{abstract} We define and characterize the notion of $(i,j)$-Baireness for bilocales. We also give internal properties of $(i,j)$-Baire bilocales which are not translated from properties of $(i,j)$-Baireness in bispaces. It turns out $(i,j)$-Baire bilocales are conservative in bilocales, in the sense that a bitopological space is almost $(i,j)$-Baire if and only if the bilocale it induces is $(i,j)$-Baire. Furthermore, in the class of Noetherian bilocales, $(i,j)$-Baireness of a bilocale coincides with $(i,j)$-Baireness of its ideal bilocale. We also consider relative versions of $(i,j)$-Baire where we show that a bilocale is $(i,j)$-Baire only if the subbilocale induced by the Booleanization is $(i,j)$-Baire. We use the characterization of $(i,j)$-Baire bilocales to introduce and characterize $(\tau_{i},\tau_{j})$-Baireness in the category of topobilocales.
	\end{abstract}
	%%%------------------------------------------------------------------
	
	\maketitle
	
	%--------------------------------------------------------------
	
	\section*{Introduction}\label{sect0}
	
In classical topology, a space is called \textit{Baire} if the intersection of every sequence of dense open sets is dense. In bispaces, an \textit{almost $(i,j)$-Baire} bispace refers to a bispace $(X,\tau_{1},\tau_{2})$ in which the intersection of any sequence of $\tau_{i}$-dense $\tau_{j}$-open subsets is $\tau_{i}$-open. A \textit{Baire} locale was introduced by Isbell \cite{I} as one in which every non-void open sublocale is of second category. To our knowledge, $(i,j)$-Baireness has not yet appeared in the category of bilocales. In this paper, we introduce and study $(i,j)$-Baire bilocales. Our definition is rather an extension of almost $(i,j)$-Baire bispaces instead of Baire locales, with the prefix ``almost" being dropped. Since the definition of almost $(i,j)$-Baire bispace is purely in terms of open subsets, we extend it to bilocales almost verbatim. We aim to extend some known bispaces results and also give some natural properties of $(i,j)$-Baire bilocales. Some of the natural results include $(i,j)$-Baireness of both the the ideal bilocale and the subbilocale induced by the smallest dense sublocale. These may lead to new results in bilocales and to strengthening some of the existing ones.

	This paper is organized as follows. Section one consists of the necessary background. In section two we introduce and characterize $(i,j)$-Baire bilocales. We also show that the class of $(i,j)$-Baire bilocales includes the following classes: (i) compact $i$-prefit bilocales, (ii) bilocales $(L,L_{1},L_{2})$ where there is an $i$-prefit compactification $h:(M,M_{1},M_{2})\rightarrow(L,L_{1},L_{2})$ with which $h_{*}[L]$ is $i$-G$_{\delta}$-dense in $M$, and (iii) $i$-pseudocomplete bilocales. In the class of Noetherian bilocales, a bilocale is $(i,j)$-Baire if and only if the induced ideal bilocale is $(i,j)$-Baire. In section three, we investigate relative versions of $(i,j)$-Baireness. We show that a bilocale is $(i,j)$-Baire only if the subbibolocale induced by the smallest dense sublocale is $(i,j)$-Baire. We also introduce and characterize relatively $(i,j)$-Baire subbilocales. It turns out that in a class of dense subbilocales, $(i,j)$-Baire coincides with relatively $(i,j)$-Baire. In section four, we define and characterize $(\tau_{i},\tau_{j})$-Baire topobilocales.
	
	\section{Preliminaries}\label{sect1}
	The book \cite{PP1} is our main reference for notions of locales and sublocales. See \cite{BBH,PP2, N2} for the theory of bilocales.
	\subsection{Locales }\label{subsect11}
	
	A \textit{locale} $L$ is a complete lattice in which
	
	$$a\wedge\bigvee B=\bigvee\{a\wedge b:b\in B\}$$ for all $a\in L$, $B\subseteq L$. 
	$1_{L}$ and $0_{L}$, with subscripts dropped if there is no possibility of confusion, respectively denote the top element and the bottom element of a locale $L$. By a \textit{point} of a locale $L$ we mean an element $a$ of $L$ such that $a\neq1$ and $b\wedge c\leq a$ implies $b\leq a$ or $c\leq a$ for all $b,c\in L$. We denote by $a^{*}$ the \textit{pseudocomplement} of an element $a\in L$. An element $a\in L$ is said to be \textit{dense} and \textit{complemented} in case $a^{*}=0$ and $a\vee a^{\ast}=1$, respectively. An element $x\in L$ is \textit{compact} if $x\leq \bigvee A$ for $A\subseteq L$ implies $x\leq \bigvee F$ for some finite $F\subseteq A$.  A locale $L$ in which the top element is compact is  called is a \textit{compact} locale. A \textit{regular} locale is a locale $L$ in which $$a=\bigvee\{x\in L:x\prec a\}$$ for every $a\in L$, where $x\prec a$ means that $x^{*}\vee a=1$.
	
	By a \textit{subframe} of a locale $L$, we mean a subset which is closed under joins and finite meets.
	
	We denote by $\mathfrak{O}X$ the locale of open subsets of a topological space $X$. 
	
	A \textit{localic map} is an infima-preserving function $f:L\rightarrow M$ between locales such that the corresponding left adjoint $f^{*}$, called the \textit{frame homomorphism}, preserves binary meets.  A frame homomorphism $h: M\rightarrow L$ is \textit{dense} if $h(x)=0$ implies $x=0$ for all $x\in M$.
	
	A \textit{sublocale} of a locale $L$ is a subset $S$ which is closed under all meets and $x\rightarrow s\in S$ for every $x\in L$ and $s\in S$, where $\rightarrow$ is a \textit{Heyting operation} on $L$ satisfying that $$a\leq b\rightarrow c \text{ if and only if } a\wedge b\leq c$$ for all $a,b,c\in L$. We denote by $\mathsf{O}$ the smallest sublocale of a locale $L$. We use $\mathcal{S}(L)$ to represent the coframe of sublocales of a locale $L$. For each $S\in \mathcal{S}(L)$, we define $$L\smallsetminus S:=\bigvee\{T\in\mathcal{S}(L):T\cap S=\mathsf{O}\}.$$
	
	The sublocales
	$${\mathfrak{c}}(a)=\{x\in L:a\leq x\}\quad\text{and}\quad \mathfrak{o}(a)=\{a\rightarrow x:x\in L\},$$ of a locale $L$ are respectively the \textit{closed} and \textit{open} sublocales induced by an element $a$ of $L$. They are complements of each other. The smallest closed sublocale of $L$ containing a sublocale $S$ is called the \textit{closure} of $S$ and denoted by $\overline{S}$ or $\cl_{L}(S)$ with subscript $L$ dropped when the locale is clear from the context. 
	
	For every $\Lambda$, $$\mathfrak{c}\left(\bigvee_{\alpha\in \Lambda}x_{\alpha}\right)=\bigcap_{\alpha\in \Lambda}\mathfrak{c}(x_{\alpha}).$$
	
	A sublocale $S$ of a locale $L$ is \textit{dense} and \textit{nowhere dense} if $\overline{S}=L$ and $S\cap\mathfrak{B}L=\mathsf{O}$, respectively, where 
	$\mathfrak{B}(L)=\{x\rightarrow 0:x\in L\}$ is the \textit{smallest dense} sublocale of $L$. By a \textit{G$_{\delta}$-sublocale} of a locale $L$, we mean a sublocale of the form $S=\bigwedge_{n\in\mathbb{N}}\mathfrak{o}(x_{n})$.
	
	For each sublocale $S\subseteq L$ there is an onto frame homomorphism $\nu_{S}:L\rightarrow S$ defined by $\nu_{S}(a)={\bigwedge}\{s\in S: a\leq s\}.$ Open sublocales and closed sublocales of a sublocale $S$ of $L$ are given by $$\mathfrak{o}_{S}(\nu_{S}(a))=S\cap \mathfrak{o}(a)\quad\text{and}\quad {\mathfrak{c}_{S}}(\nu_{S}(a))=S\cap \mathfrak{c}(a),$$ respectively, for $a\in L$. %For any sublocale $S$ of a locale $L$ and $x\in L$, $S\subseteq\mathfrak{o}(x)$ if and only if $\nu_{S}(x)=1$.
	
	Each localic map $f:L\rightarrow M$ induces the functions $f[-]:\mathcal{S}(L)\rightarrow \mathcal{S}(M)$ given by the set-theoretic image of each sublocale of $L$ under $f$, and   $f_{-1}[-]:\mathcal{S}(M)\rightarrow\mathcal{S}(L)$ given by $$f_{-1}[T]={\bigvee}\{A\in \mathcal{S}(L):A\subseteq f^{-1}(T)\}.$$ For a localic map $f:L\rightarrow M$ and $x\in M$,  $$f_{-1}[\mathfrak{c}_{M}(x)]=\mathfrak{c}_{L}(h(x))\quad\text{and}\quad f_{-1}[\mathfrak{o}_{M}(x)]=\mathfrak{o}_{L}(h(x)).$$
	
	We denote by $\widetilde{A}$ the sublocale of $\mathfrak{O}X$ \textit{induced} by a subset $A$ of a topological space $X$. 
	
	%###############
	%###############
	\subsection{Bilocales}
	
	A \textit{bilocale} is a triple $(L,L_{1}, L_{2})$ where $L_{1},L_{2}$ are subframes of a locale $L$ and for all $a\in L$, \[ a=\bigvee\{a_{1}\wedge a_{2}: a_{1}\in L_{1}, a_{2}\in L_{2}\text{ and } a_{1}\wedge a_{2}\leq a\}.\] 
	
	We call $L$ the \textit{total part} of $(L,L_{1}, L_{2})$, and $L_{1}$ and $L_{2}$ the first and second parts, respectively. We use the notations $L_{i},L_{j}$ to denote the first or second parts of $(L,L_{1},L_{2})$, always assuming that $i,j=1,2$, $i\neq j$. 
	
	For every bispace $(X,\tau_{1},\tau_{2})$, there is a corresponding bilocale $(\tau_{1}\vee\tau_{2},\tau_{1},\tau_{2})$.
	
	The \textit{bilocale pseudocomplement} of $c\in L_{i}$ is given by $$c^{\bullet}=\bigvee\{x\in L_{j}:x\wedge c=0\}.$$
	For all $a\in L_{j},b\in L_{i}$, $a\wedge b=0$ if and only if $a\leq b^{\bullet}$.
	
	A bilocale $(L,L_{1},L_{2})$ is \textit{compact} if its total part is compact, and \textit{regular} provided that $$x=\bigvee\{a\in L_{i}:a\prec_{i} x\}$$ for every $x\in L_{i}$, where $a\prec_{i} x$ means that there is $c\in L_{j}$ such that $a\wedge c=0$ and $c\vee x=1$.

	A \textit{biframe homomorphism} (or \textit{biframe map}) $h:(M,M_{1},M_{2})\rightarrow(L,L_{1},L_{2})$ is a frame homomorphism $h:M\rightarrow L$ for which $h(M_{i})\subseteq L_{i}\quad (i=1,2).$ A biframe map $h:(M,M_{1},M_{2})\rightarrow(L,L_{1},L_{2})$ is \textit{dense} if $h: M\rightarrow L$ is dense. It is \textit{onto} if $h[M_{i}]=L_{i}$ for $i=1,2$.
	
	A \textit{subbilocale} of a bilocale $(L,L_{1},L_{2})$ is a triple $(S,S_{1},S_{2})$ where $S$ is a sublocale of $L$ and $$S_{i}=\nu_{S}[L_{i}]\quad \text{for}\quad i=1,2.$$ We shall say that $(S,S_{1},S_{2})$ is a $P$-subbilocale in case $S$ has property $P$.
	
	Recall that for a bilocale $(L,L_{1},L_{2})$ and a sublocale $S$ of $L$:  \cite{N2}
	$$\Int_{i}(S)=\bigvee\{\mathfrak{o}(a):a\in L_{i}, \mathfrak{o}(a)\subseteq S\}\quad(i=1,2).$$
	and \cite{PP2} $$\cl_{i}(S)=\bigcap\{\mathfrak{c}(a):a\in L_{i}, S\subseteq\mathfrak{c}(a)\}=\mathfrak{c}\left(\bigvee\{a\in L_{i}:S\subseteq\mathfrak{c}(a)\}\right)\quad(i=1,2).$$
	
For each $a\in L_{i}$, $\mathfrak{c}(a^{\bullet})=\cl_{j}(\mathfrak{o}(a))$.

		A sublocale $A$ of a bilocale $(L,L_{1},L_{2})$ is \textit{$i$-dense} if $\cl_{i}(A)=L$. This is equivalent to saying that $S$ is $i$-dense if and only if for each non-zero $x\in L_{i}$, $\mathfrak{o}(x)\cap S\neq \mathsf{O}$. Every sublocale containing an $i$-dense sublocale is $i$-dense. 
		
		Given a bilocale $(L,L_{1},L_{2})$, a sublocale $S$ of $L$ is \textit{$(i,j)$-nowhere dense} if $\Int_{j}(\cl_{i}(S))=\mathsf{O}$ $(i\neq j\in\{1,2\})$. As a result, a sublocale $S$ of a bilocale $(L,L_{1},L_{2})$ is $(i,j)$-nowhere dense if and only if  $L\smallsetminus\cl_{i}(S)$ is $j$-dense if and only if $\overline{S}$ is $(i,j)$-nowhere dense.  
		Furthermore, an element $a\in L_{i}$ is $j$-dense if and only if $\mathfrak{c}(a)$ is $(i,j)$-nowhere dense.
		
		When dealing with subbilocales, say $(S,S_{1},S_{2})$, we write $i_{S}$-dense, $i_{S}$-open and $(i_{S},j_{S})$-nowhere dense instead of $i$-dense, $i$-open and $(i,j)$-nowhere dense. 
		
		By an \textit{$i$-G$_{\delta}$-sublocale} of a bilocale $(L,L_{1},L_{2})$, we mean a sublocale of the form $S=\bigwedge_{n\in\mathbb{N}}\mathfrak{o}(x_{n})$ where each $x_{n}\in L_{i}$. A sublocale of a bilocale $(L,L_{1},L_{2})$ is \textit{$i$-$G_{\delta}$-dense} if it meets every nonvoid $i$-$G_{\delta}$-sublocales.
		
		For the bilocale $(\tau_{1}\vee\tau_{2},\tau_{1},\tau_{2})$ induced by a bispace $(X,\tau_{1},\tau_{2})$, we shall write $U$ is $\tau_{i}$-open if $U\in \tau_{i}$ and $U$ is $\tau_{i}$-dense if $U$ is dense with respect to the topological space $(X,\tau_{i})$.
		
		%###############
	%###############
	
\section{$(i,j)$-Baire Bilocales}

Recall that a bispace $(X,\tau_{1},\tau_{2})$ is \textit{almost $(i,j)$-Baire} \cite{D} if any collection $\{U_{n}:n\in \mathbb{N}\}$ of $\tau_{i}$-dense $\tau_{j}$-open subsets of $X$ satisfies the condition $\bigcap_{n\in\mathbb{N}}U_{n}$ is $\tau_{i}$-dense. In this section, we extend the definition of almost $(i,j)$-Baire bispaces to bilocales where the prefix ``almost" shall be dropped. We aim to define $(i,j)$-Baire bilocales in such a way that a bispace $(X,\tau_{1},\tau_{2})$ is almost $(i,j)$-Baire if and only if the bilocale $(\tau_{1}\vee\tau_{2},\tau_{1},\tau_{2})$ is $(i,j)$-Baire. 

We shall call an open (resp. closed) sublocale \textit{$i$-open} (resp. $i$-closed) in case the inducing element is an element of $L_{i}$. 
	
	\begin{definition}
		A bilocale $(L,L_{1},L_{2})$ is said to be \textit{$(i,j)$-Baire} if the intersection of countably many $i$-dense $j$-open sublocales is $i$-dense. 
	\end{definition}

\begin{example}
	Since, in locales, the intersection of dense sublocales is dense, every symmetric bilocale (bilocale of the form $(L,L,L)$) is $(i,j)$-Baire.
\end{example}

For a bispace $(X,\tau_{1},\tau_{2})$ and $A\subseteq X$, define $$\widetilde{A}=\{\inte_{\tau_{1}\vee\tau_{2}}((X\smallsetminus A)\cup G):G\in\tau_{1}\vee\tau_{2}\}.$$ It is clear that $\widetilde{A}$ is a sublocale of $\tau_{1}\vee\tau_{2}$. For each $x\in X$, $\widetilde{x}=X\smallsetminus\cl_{\tau_{1}\vee\tau_{2}}{\{x\}}$ is a point of $\tau_{1}\vee\tau_{2}$. Just like in the case of locales, $\mathfrak{o}(U)=\widetilde{U}$ for every $U\in\tau_{1}\vee\tau_{2}$. 

Recall from \cite{L} that given a topological property $P$, a bispace $(X,\tau_{1},\tau_{2})$ is sup-$P$ if $(X,\tau_{1}\vee\tau_{2})$ has property $P$. In \cite{N1}, we proved the following result.
\begin{lemma}
	 Let $(X,\tau_{1},\tau_{2})$ be a sup-T$_{D}$-bispace. Then $A\subseteq X$ is $\tau_{i}$-dense in $(X,\tau_{1},\tau_{2})$ iff $\widetilde{A}$ is $i$-dense in $(\tau_{1}\vee\tau_{2},\tau_{1},\tau_{2})$. 
\end{lemma}

In \cite{PP0}, the authors show that if $X$ is a topological space and $Y\subseteq X$, then $$\widetilde{Y}=\bigvee\{\{X\smallsetminus\overline{\{y\}},1_{\mathfrak{O}X}\}:y\in Y\}.$$
\begin{proposition}
	Let $(X,\tau_{1},\tau_{2})$ be a sup-T$_{D}$-bispace in which G$_{\delta}$-sublocales of $(\tau_{1}\vee \tau_{2},\tau_{1},\tau_{2})$ are complemented. Then $(X,\tau_{1},\tau_{2})$ is almost $(i,j)$-Baire iff $(\tau_{1}\vee \tau_{2},\tau_{1},\tau_{2})$ is $(i,j)$-Baire.
\end{proposition}
\begin{proof}
	$(\Longrightarrow)$: Let $\{\mathfrak{o}(U_{n}):n\in\mathbb{N}\}$ be a collection of $i$-dense $j$-open sublocales. It follows that $\{U_{n}:n\in\mathbb{N}\}$ is a collection of $\tau_{i}$-dense $\tau_{j}$-open subsets of $X$. It follows that $\bigcap_{n\in\mathbb{N}}U_{n}$ is $\tau_{i}$-dense.
	
	Claim: $\bigwedge_{n\in\mathbb{N}}\mathfrak{o}(U_{n})$ is $i$-dense.
	
	Proof: Let $\mathfrak{o}(V)$ be an $i$-open sublocale such that $\mathfrak{o}(V)\cap\bigwedge_{n\in\mathbb{N}}\mathfrak{o}(U_{n})=\mathsf{O}$. Then $$\bigwedge_{n\in\mathbb{N}}\mathfrak{o}(V\cap U_{n})=\mathsf{O}.$$ We must have that $\bigcap_{n\in\mathbb{N}}(V\cap U_{n})=\emptyset$. Otherwise, there is $x\in V\cap U_{n}$ for each $n\in\mathbb{N}$. Since each $V\cap U_{n}$ is $\tau_{1}\vee\tau_{2}$, $$\widetilde{V\cap U_{n}}=\mathfrak{o}(V\cap U_{n})\ni\widetilde{x}$$ for each $n\in\mathbb{N}$. Therefore $\widetilde{x}\in \bigwedge_{n\in\mathbb{N}}\mathfrak{o}(V\cap U_{n})$ which is impossible. Therefore $V=\emptyset$ so that $\mathfrak{o}(V)=\mathsf{O}$. Thus  $\bigwedge_{n\in\mathbb{N}}\mathfrak{o}(U_{n})$ is $i$-dense.
	
	$(\Longleftarrow)$: Let $\{U_{n}:n\in\mathbb{N}\}$ be a collection of $\tau_{i}$-dense $\tau_{j}$-open subsets of $X$. Then $\{\mathfrak{o}(U_{n}):n\in\mathbb{N}\}$ is a collection of $i$-dense $j$-open sublocales of $\tau_{1}\vee\tau_{2}$. It follows that $\bigwedge_{n\in\mathbb{N}}\mathfrak{o}(U_{n})$ is $\tau_{i}$-dense. To show that $\bigcap_{n\in\mathbb{N}}U_{n}$ is $\tau_{i}$-dense, let $V$ be a nonempty $\tau_{i}$-open subset of $X$ such that $V\cap\left(\bigcap_{n\in\mathbb{N}}U_{n}\right)=\emptyset$. Then $$\bigcup_{n\in\mathbb{N}} (X\smallsetminus (V\cap U_{n}))=X.$$
	
	Observe that $\bigvee_{n\in\mathbb{N}}\mathfrak{c}(V\cap U_{n})=\mathfrak{O}X$: Let $p\in X$. Then $p\in X\smallsetminus(V\cap U_{n})$ for some $n\in\mathbb{N}$. Therefore $\overline{\{p\}}\subseteq X\smallsetminus (V\cap U_{n})$ for some $n\in\mathbb{N}$ so that $V\cap U_{n}\subseteq X\smallsetminus\overline{\{p\}}$. This implies that $X\smallsetminus\overline{\{p\}}\in\mathfrak{c}(V\cap U_{n})$. As a result, $$\{X\smallsetminus\overline{\{p\}},1_{\tau_{1}\vee\tau_{2}}\}\subseteq\mathfrak{c}(V\cap U_{n}).$$
	
	Therefore \begin{align*}
		\tau_{1}\vee\tau_{2}&=\bigvee\{\{X\smallsetminus\overline{\{p\}},1_{\tau_{1}\vee\tau_{2}}\}:p\in X\}\\
		&\subseteq\bigvee\{\mathfrak{c}(V\cap U_{n}):{n\in\mathbb{N}}\}\\
		&\subseteq \tau_{1}\vee\tau_{2}.
	\end{align*}
Since $\bigwedge_{n\in\mathbb{N}}\mathfrak{o}(U_{n})$ is $i$-dense and $\mathfrak{o}(V)$ is non-void $i$-open, $$\mathfrak{o}(V)\cap\bigwedge_{n\in\mathbb{N}}\mathfrak{o}(U_{n})\neq\mathsf{O},$$ so that $$\mathsf{O}\neq\bigwedge_{n\in\mathbb{N}}\mathfrak{o}(V\cap U_{n}).$$ 
Because $\bigwedge_{n\in\mathbb{N}}\mathfrak{o}(V\cap U_{n})$ is a complemented G$_{\delta}$-sublocale, we have that $$\bigvee_{n\in\mathbb{N}}\mathfrak{c}(V\cap U_{n})\neq\tau_{1}\vee\tau_{2},$$ which is impossible.
\end{proof}

\begin{definition}
	Let $(L,L_{1},L_{2})$ be a bilocale. A sublocale $S$ of $L$ is said to be of  \textit{$(i,j)$-first category} if there are countably many $(i,j)$-nowhere dense sublocales $S_{n}$, $n\in\mathbb{N}$, such that $S\subseteq \bigvee_{n\in\mathbb{N}}S_{n}$. It is of \textit{ $(i,j)$-second category} if it is not of  $(i,j)$-first category.
\end{definition}
\begin{theorem}\label{mainresult}
	Let $(L,L_{1},L_{2})$ be a bilocale whose $j$-G$_{\delta}$-sublocales are complemented in $L$. The following statements are equivalent:
	\begin{enumerate}
		\item $(L,L_{1},L_{2})$ is $(i,j)$-Baire.
		\item Each non-void $i$-open sublocale is of $(j,i)$-second category.
		\item Every sublocale of $(j,i)$-first category has void $i$-interior.
		\item The supplement of every sublocale of $(j,i)$-first category is $i$-dense. 
	\end{enumerate}
\end{theorem}
\begin{proof}
	$(1)\Longrightarrow(2)$: Let $U$ be a non-void $i$-open sublocale of $L$ and assume that $U\subseteq \bigvee_{n\in\mathbb{N}}S_{n}$ for some collection $\{S_{n}:n\in\mathbb{N}\}$ of $(j,i)$-nowhere dense sublocales. Then $U\subseteq \bigvee_{n\in\mathbb{N}}\overline{S_{n}}$ where members of the collection $\{\overline{S_{n}}:n\in\mathbb{N}\}$ are $(j,i)$-nowhere dense. It follows that  members of the collection $\{L\smallsetminus \cl_{j}(\overline{S_{n}}):n\in\mathbb{N}\}$ are $i$-dense $j$-open sublocales. By hypothesis, $\bigwedge_{n\in\mathbb{N}}(L\smallsetminus\cl_{j}(\overline{S_{n}})$ is $i$-dense so that $$U\cap \left(\bigwedge_{n\in\mathbb{N}}(L\smallsetminus \cl_{j}(\overline{S_{n}}))\right)\neq\mathsf{O}.$$ Therefore \begin{align*}
		\mathsf{O}&\neq \left(\bigvee_{k\in\mathbb{N}}\overline{S_{k}}\right)\cap\left(\bigwedge_{n\in\mathbb{N}}(L\smallsetminus \cl_{j}(\overline{S_{n}}))\right)\\
		&=\bigvee_{k\in\mathbb{N}}\left(\overline{S_{k}}\cap \left(\bigwedge_{n\in\mathbb{N}}(L\smallsetminus \overline{S_{n}})\right)\right)\\
		&\subseteq\bigvee_{k\in\mathbb{N}}\left(\overline{S_{k}}\cap(L\smallsetminus\overline{S_{k}})\right)\\
		&=\mathsf{O}
	\end{align*}which is impossible.

	$(2)\Longrightarrow(3)$: Let $S$ be a sublocale of $(j,i)$-first category with a non-void $i$-interior. We then get that $\inte_{i}(S)$ is a non-void $i$-open sublocale which must be of $(j,i)$-second category by (2). This is a contradiction.
	
	$(3)\Longrightarrow(4)$: Let $S$ be a sublocale of $L$ which is of $(j,i)$-first category and choose $x\in L_{i}$ with $\mathfrak{o}(x)\cap (L\smallsetminus S)=\mathsf{O}$. Then $\mathfrak{o}(x)\subseteq S$. Since $S$ satisfies the conditions hypothesized in $(3)$, $\inte_{i}(S)\neq\mathsf{O}$, so that the $i$-open sublocale $\mathfrak{o}(x)$ is void.
	
	$(4)\Longrightarrow(1)$: Let $\{\mathfrak{o}(x_{n}):n\in\mathbb{N}\}$ be a collection of $i$-dense $j$-open sublocales and assume that there is an $i$-open sublocale $\mathfrak{o}(y)$ with $$\mathfrak{o}(y)\cap\left(\bigwedge_{n\in\mathbb{N}}\mathfrak{o}(x_{n})\right)=\mathsf{O}$$ Then $\mathfrak{o}(y)\subseteq\bigvee_{n\in\mathbb{N}}\mathfrak{c}(x_{n})$, making $\mathfrak{o}(y)$ a sublocale of $(j,i)$-first category. By (4), $L\smallsetminus \mathfrak{o}(y)=\mathfrak{c}(y)$ is $i$-dense, i.e., $\cl_{i}(\mathfrak{c}(y))=\mathfrak{c}(y)=1$. Therefore $y=0$ so that $\mathfrak{o}(y)=\mathsf{O}$. Hence $\left(\bigwedge_{n\in\mathbb{N}}\mathfrak{o}(x_{n})\right)$ is $i$-dense.
\end{proof}
We shall say that a collection $\mathcal{C}$ of sublocales of $L$ has the \textit{Finite Intersection Property (FIP)} if the intersection of every finite subcollection of $\mathcal{C}$ has a non-void intersection.

\begin{proposition}
	If a bilocale $(L,L_{1},L_{2})$ is compact, then every collection of closed sublocales with the FIP has a non-void intersection.
\end{proposition}
\begin{proof}
	Let $\{\mathfrak{c}(x_{\alpha}):\alpha\in \Lambda\}$ be a collection with the FIP and assume that $\bigcap_{\alpha\in \Lambda}\mathfrak{c}(x_{\alpha})=\mathsf{O}$. Then $$L=L\smallsetminus\bigcap_{\alpha\in \Lambda}\mathfrak{c}(x_{\alpha})=\bigvee_{\alpha\in \Lambda}\mathfrak{o}(x_{\alpha}),$$ making $\mathfrak{o}\left(\bigvee_{\alpha\in \Lambda}x_{\alpha}\right)=L$. Therefore $\bigvee_{\alpha\in \Lambda}x_{\alpha}=1$. Since $(L,L_{1},L_{2})$ is compact, there is a finite set $F\subseteq \Lambda$ such that $\bigvee_{\alpha\in F}x_{\alpha}=1$. We get that $$\mathsf{O}=\mathfrak{c}\left(\bigvee_{\alpha\in F}x_{\alpha}\right)=\bigcap_{\alpha\in F}\mathfrak{c}(x_{\alpha}),$$ which contradicts that $\{\mathfrak{c}(x_{\alpha}):\alpha\in \Lambda\}$ has the FIP.
\end{proof}
\begin{observation}
	The converse of the preceding result holds. We are however interested in the forward direction, that is why we only proved it.
\end{observation}
Recall from \cite{PP} that a locale $L$ is \textit{prefit} if for each nonzero $x\in L$ there is a nonzero $y\in L$ such that $y^{\star}\vee x=1$. A bispace $(X,\tau_{1},\tau_{2})$ is \textit{almost regular} if for each nonempty $U\in \tau_{i}$, there is nonempty $V\in\tau_{i}$ such that $\cl_{\tau_{j}}(V)\subseteq U$. Since prefitness is a localic version of almost regularity in spaces (spaces in which every nonempty open set contains some closure of a nonempty open subset), we define a \textit{prefit bilocale} $(L,L_{1},L_{2})$ using the notion of almost regular bispace as one in which for each $x\in L_{i}$, $i=1,2$, there is $y\in L_{i}$ such that $y^{\bullet}\vee x=1$. Related to prefit bilocales, we give the following definition.
\begin{definition}
	Call a bilocale $(L,L_{1},L_{2})$ \textit{strongly prefit} if for every nonzero $x\in L$ there are nonzero $a\in L_{i}$ and $b\in L_{j}$ such that $a^{\bullet}\vee x=1=b^{\bullet}\vee x$. We shall say that $(L,L_{1},L_{2})$ is \textit{$i$-prefit} in case for every nonzero $x\in L$, there is a nonzero $y\in L_{i}$ such that $y^{\bullet}\vee x=1$.
\end{definition}

We consider some examples.

\begin{example}
	\begin{enumerate}
		\item Every strongly prefit bilocale is both $i$-prefit and prefit. 
		\item The biframe of reals is an example of a prefit bilocale which is not $i$-prefit and hence not strongly prefit by (1). 
		\item For any almost regular bispace $(X,\tau_{1},\tau_{2})$ with $\tau_{1}\subseteq\tau_{2}$, the bilocale $(\tau_{1}\vee\tau_{2},\tau_{1},\tau_{2})$ is $2$-prefit. In particular, if $L$ is prefit, then $(L,L,L)$ is strongly prefit and hence $i$-prefit ($i=1,2$).
		\item By \cite{S}, a bilocale $(L,L_{1},L_{2})$ is Boolean if for each $x\in L_{i}$, $i=1,2$, there is $c\in L_{j}$ $(i\neq j)$ such that $x\wedge c=0$ and $x\vee c=1$. Boolean and $i$-prefit are incomparable: Consider the set $X=\{a,b,c,d\}$ endowed with topologies $\tau_{1}=\{\emptyset,X,\{a\},\{b\},\{a,b\}\}$ and $\tau_{2}=\{\emptyset,X,\{b,c,d\},\{a,c,d\},\{c,d\}\}$. It is clear that $(\tau_{1}\vee\tau_{2},\tau_{1},\tau_{2})$ is Boolean. This bilocale is not $i$-prefit ($i=1,2$) since for the set $\{a\}\in\tau_{1}\vee\tau_{2}$, there is no nonempty $U\in\tau_{i}$ satisfying that $U^{\bullet}\vee \{a\}=X$.
		
		For any non-Boolean prefit locale $L$, $(L,L,L)$ is an example of a non-Boolean strongly prefit and hence $i$-prefit ($i=1,2$) bilocale.
	\end{enumerate}
\end{example}
In the following result, we show that the class of $(i,j)$-Baire bilocales contains compact $i$-prefit bilocales.
\begin{proposition}\label{jprefitbaire}
	Every compact $i$-prefit bilocale is $(i,j)$-Baire.
\end{proposition}
\begin{proof}
	Let $(L,L_{1},L_{2})$ be a compact $i$-prefit bilocale and choose a collection $\{\mathfrak{o}(x_{n}):n\in\mathbb{N}, x_{n}\in L_{i}\}$ of $i$-dense $j$-open sublocales and a non-void $j$-open sublocale $\mathfrak{o}(y)$. Then $$\mathfrak{o}(y)\cap\mathfrak{o}(x_{n})\neq\mathsf{O}$$ for each $n\in \mathbb{N}$. This makes $y\wedge x_{n}\neq 0$. Since $(L,L_{1},L_{2})$ is $i$-prefit, there is nozero $b_{1}\in L_{i}$ such that $b_{1}^{\bullet}\vee (y\wedge x_{1})=1$. Because $\mathfrak{o}(x_{2})$ is $i$-dense, we have that $\mathfrak{o}(x_{2})\cap \mathfrak{o}(b_{1})\neq\mathsf{O}$ so that $x_{2}\wedge b_{1}$ is a nonzero element of $L$. By $i$-prefitness again, there is nonzero $b_{2}\in L_{i}$ such that $b_{2}^{\bullet}\vee (x_{2}\wedge b_{1})=1$. Continuing like this for $n=3,4,..$, we find $b_{n}\in L_{i}$ such that $b_{n}^{\bullet}\vee (x_{n}\wedge b_{n-1})=1$. Therefore $\mathfrak{c}(b_{n}^{\bullet})\subseteq \mathfrak{o}(x_{n})\cap\mathfrak{o}(b_{n-1})$. Since each $b_{n}\in L_{i}$, we have that $\mathfrak{c}(b_{n}^{\bullet})=\cl_{j}(\mathfrak{o}(b_{n}))$. Therefore $$...\subseteq\mathfrak{c}(b_{3}^{\bullet})=\cl_{j}(\mathfrak{o}(b_{3}))\subseteq\mathfrak{c}(b_{2}^{\bullet})=\cl_{j}(\mathfrak{o}(b_{2}))\subseteq \mathfrak{c}(b_{1}^{\bullet})=\cl_{j}(\mathfrak{o}(b_{1}))\subseteq \mathfrak{o}(y)\cap\mathfrak{o}(x_{1}).$$ We now have the decreasing sequence $\mathfrak{c}(b_{1}^{\bullet}), \mathfrak{c}(b_{2}^{\bullet}), \mathfrak{c}(b_{3}^{\bullet}),...$ of closed sublocales, so that the collection $\{\mathfrak{c}(b_{n}^{\bullet}):n\in\mathbb{N}\}$ has the FIP. By compactness of $(L,L_{1},L_{2})$, $\bigcap_{n\in\mathbb{N}}\mathfrak{c}(b_{n}^{\bullet})\neq\mathsf{O}$. Because $$\bigcap_{n\in\mathbb{N}}\mathfrak{c}(b_{n}^{\bullet})\subseteq \mathfrak{o}(y)\cap\bigcap_{n\in\mathbb{N}}\mathfrak{o}(x_{n}),$$ we then have that $$\mathfrak{o}(y)\cap\bigcap_{n\in\mathbb{N}}\mathfrak{o}(x_{n})\neq\mathsf{O},$$ making $\bigcap_{n\in\mathbb{N}}\mathfrak{o}(x_{n})$ an $i$-dense sublocale.
\end{proof}
The converse of Proposition \ref{jprefitbaire} is not always true, as shown below.
\begin{example}
	 For the bispace $(\mathbb{N},\tau_{D},\tau_{cf})$,  $(\tau_{D}\vee\tau_{cf}=\tau_{D},\tau_{D},\tau_{cf})$ is $(\tau_{D},\tau_{cf})$-Baire but not compact. Its $(\tau_{D},\tau_{cf})$-Baireness follows since $\mathbb{N}$ is the only $\tau_{D}$-dense member of $\tau_{cf}$.
\end{example}
Recall that for any onto frame homomorphism $h:M\rightarrow L$, $h_{*}:L\rightarrow h_{*}[L]$ is a frame isomorphism. If $h$ is further dense, then $h_{*}(0)=0$.
\begin{lemma}\label{frameiso}
	Let $h:(M,M_{1},M_{2})\rightarrow(L,L_{1},L_{2})$ be a dense and onto biframe map. If $x\in L_{i}$ is $j$-dense, then $h_{*}(x)\wedge a\neq 0$ whenever $a\in M_{i}$ is nonzero. 
\end{lemma}
\begin{proof}
	Let $a\in M_{i}$ be nonzero such that $a\wedge h_{*}(x)=0$. Then $$0=h(a)\wedge h(h_{*}(x))=h(a)\wedge x.$$ Since $h[M_{i}]\subseteq L_{i}$, $h(a)\in L_{i}$ so that $h(a)=0$. Therefore $a\leq h_{*}(h(a))=h_{*}(0)=0$. 
\end{proof}
\begin{proposition}
	Let $(L,L_{1},L_{2})$ be a bilocale. If there is a dense onto biframe map $h:(M,M_{1},M_{2})\rightarrow(L,L_{1},L_{2})$ from an $(i,j)$-Baire bilocale $(M,M_{1},M_{2})$ with which $h_{*}[L]$ is $i$-G$_{\delta}$-dense in $M$, then $(L,L_{1},L_{2})$ is $(i,j)$-Baire.
\end{proposition}
\begin{proof}
	Suppose that the hypothesized statement is true and let $\{\mathfrak{o}(x_{n}):n\in\mathbb{N}\}$ be a collection of $i$-dense $j$-open sublocales of $L$.  Now, if $\mathfrak{o}(y)\cap\left(\bigwedge_{n\in\mathbb{N}}\mathfrak{o}(x_{n})\right)=\mathsf{O}$ for some $i$-open sublocale $\mathfrak{o}(y)$ of $L$, then $$\mathsf{O}=h_{*}[\mathfrak{o}(y)]\cap h_{*}\left[\bigwedge_{n\in\mathbb{N}}\mathfrak{o}(x_{n})\right]=h_{*}[\mathfrak{o}(y)]\cap \left(\bigwedge_{n\in\mathbb{N}}h_{*}[\mathfrak{o}(x_{n})]\right)$$ where the first equality follows since $h_{*}[\mathsf{O}]=\mathsf{O}$ and $h_{*}$ is injective, and the second equality follows since the total part $h_{*}:L\rightarrow M$ is a right adjoint. By virtue of $h_{*}:L\rightarrow h_{*}[L]$ being a frame isomorphism and hence open, we get that $$\mathfrak{o}_{h_{*}[L]}(h_{*}(y))\cap \left(\bigwedge_{n\in\mathbb{N}}\mathfrak{o}_{h_{*}[L]}(h_{*}(x_{n}))\right)=\mathsf{O}.$$  For each $n\in\mathbb{N}$, $h_{*}(x_{n})=h_{*}(h(a_{n}))$ for some $a_{n}\in M_{i}$. Therefore \begin{align*}
		\mathsf{O}&=\mathfrak{o}_{h_{*}[L]}(h_{*}(y))\cap \left(\bigwedge_{n\in\mathbb{N}}\mathfrak{o}_{h_{*}[L]}(h_{*}(h(a_{n})))\right)\\
		&={h_{*}[L]}\cap \mathfrak{o}(h_{*}(y))\cap \left(\bigwedge_{n\in\mathbb{N}}\left({h_{*}[L]}\cap\mathfrak{o}(h_{*}(h(a_{n})))\right)\right)\\
		&=h_{*}[L]\cap \mathfrak{o}(h_{*}(y))\cap \left(\bigwedge_{n\in\mathbb{N}}\mathfrak{o}(h_{*}(h(a_{n})))\right)\\
		&\supseteq h_{*}[L]\cap \mathfrak{o}(h_{*}(y))\cap \left(\bigwedge_{n\in\mathbb{N}}\mathfrak{o}(a_{n})\right).
	\end{align*}  Since $i$-G$_{\delta}$-dense sublocales are dense, $h_{*}[L]$ is a dense sublocale of $M$. Each of the $a_{n}$'s is $i$-dense: Pick $b\in M_{i}$ such that $b\wedge a_{n}=0$. Then $h(b)\wedge h(a_{n})=0$ so that $$h_{*}(h(b))\wedge h_{*}(h(a_{n}))=h_{*}(0)=0,$$ where the latter equality follows since $h_{*}:L\rightarrow M$ is a dense localic map. Therefore $b\wedge h_{*}(x_{n})=0$. By Lemma \ref{frameiso}, $b=0$ and hence each $a_{n}$ is $i$-dense. 

Therefore the collection $\{\mathfrak{o}(a_{n}):n\in\mathbb{N}\}$ consists of $i$-dense $j$-open sublocales of $M$. We then get that $\mathfrak{o}(h_{*}(y))\cap \left(\bigwedge_{n\in\mathbb{N}}\left(\mathfrak{o}(a_{n})\right)\right)$ is an $i$-G$_{\delta}$-sublocale. Because $h_{*}[L]$ is $i$-G$_{\delta}$-dense, so $$\mathfrak{o}(h_{*}(y))\cap \left(\bigwedge_{n\in\mathbb{N}}\mathfrak{o}(a_{n})\right)=\mathsf{O}.$$ Since $(M,M_{1},M_{2})$ is $(i,j)$-Baire, it follows that $\bigwedge_{n\in\mathbb{N}}\mathfrak{o}(a_{n})$ is $i$-dense, so that $\mathfrak{o}(h_{*}(y))=\mathsf{O}$.  Therefore $h_{*}(y)=0$, so that $$0=h(h_{*}(y))=y.$$This means that $\mathfrak{o}(y)=\mathsf{O}$. Thus $\bigwedge_{n\in\mathbb{N}}\mathfrak{o}(x_{n})$ is $i$-dense, and hence $(L,L_{1},L_{2})$ is $(i,j)$-Baire.
\end{proof}
A \textit{compactification} of a bilocale $(L,L_{1},L_{2})$ is a dense and onto biframe map $h:(M,M_{1},M_{2})\rightarrow(L,L_{1},L_{2})$ from a compact regular bilocale $(M,M_{1},M_{2})$. So, for a bilocalic property $P$, we shall say that $(L,L_{1},L_{2})$ has a $P$-compactification in case $(M,M_{1},M_{2})$ has property $P$.
\begin{corollary}
	Let $(L,L_{1},L_{2})$ be a bilocale. If there is a $i$-prefit compactification $h:(M,M_{1},M_{2})\rightarrow(L,L_{1},L_{2})$ with which $h_{*}[L]$ is $i$-G$_{\delta}$-dense in $M$, then $(L,L_{1},L_{2})$ is $(i,j)$-Baire.
\end{corollary}
\begin{definition}
	Let $(L,L_{1},L_{2})$ be a bilocale. An \textit{$i$-$\pi$-base} for $(L,L_{1},L_{2})$ is a collection $\mathcal{C}$ of non-void $i$-open sublocales such that each non-void $i$-open sublocale of $L$ contains a member of $\mathcal{C}$. A bilocale is said to be \textit{$i$-pseudocomplete} if it is $i$-prefit and it has a sequence $(\mathcal{C}_{n})_{n\in\mathbb{N}}$ of $i$-$\pi$-bases  such that whenever $\mathfrak{o}(x_{n})\in\mathcal{C}_{n}$ and $\cl_{j}(\mathfrak{o}(x_{n+1}))\subseteq\mathfrak{o}(x_{n})$ for each $n$, then $\bigcap_{n\in\mathbb{N}}\mathfrak{o}(x_{n})\neq\mathsf{O}$.
\end{definition}
\begin{proposition}
	Every $i$-pseudocomplete bilocale is $(i,j)$-Baire.
\end{proposition}
\begin{proof}
	Let $(L,L_{1},L_{2})$ be a pseudocomplete bilocale and pick a collection $\{\mathfrak{o}(x_{n}):n\in\mathbb{N}\}$ of $i$-dense $j$-open sublocales. Since $(L,L_{1},L_{2})$ is pseudocomplete, there is a sequence $(\mathcal{C}_{n})_{n\in\mathbb{N}}$ of $i$-$\pi$-bases with the corresponding pseudocompleteness property. For each non-void $i$-open sublocale $\mathfrak{o}(y)$, we have that each $\mathfrak{o}(y)\cap \mathfrak{o}(x_{n})$ is a non-void open sublocale, so that $y\wedge x_{n}\neq0$. Since $(L,L_{1},L_{2})$ is $i$-prefit, for $n=1$, there is nonzero $a_{1}\in L_{i}$ such that $a_{1}^{\bullet}\vee (y\wedge x_{1})=1$. This makes $$\mathsf{O}\neq\mathfrak{o}(a_{1})\subseteq\cl_{j}(\mathfrak{o}(a_{1}))\subseteq \mathfrak{o}(y)\cap\mathfrak{o}(x_{1}).$$ 
	For the $i$-$\pi$-base $\mathcal{C}_{1}$, there is $\mathfrak{o}(c_{1})\in \mathcal{C}_{1}$ with $$\mathfrak{o}(c_{1})\subseteq\cl_{j}(\mathfrak{o}(c_{1}))\subseteq \mathfrak{o}(a_{1})\subseteq\mathfrak{o}(y)\cap\mathfrak{o}(x_{1}).$$ Using the fact that $\mathfrak{o}(x_{2})$ is $i$-dense and $\mathfrak{o}(c_{1})$ is non-void $i$-open, we get that $\mathfrak{o}(c_{1})\cap\mathfrak{o}(x_{2})$ is a non-void open sublocale. Since $(L,L_{1},L_{2})$ is $i$-prefit, there is nonzero $a_{2}\in L_{i}$ such that $$\mathsf{O}\neq\mathfrak{o}(a_{2})\subseteq\cl_{j}(\mathfrak{o}(a_{2}))\subseteq \mathfrak{o}(c_{1})\cap\mathfrak{o}(x_{2}).$$ An application of pseudocompleteness to the $i$-$\pi$-base $\mathcal{C}_{2}$ yields an existence of $\mathfrak{o}(c_{2})\in \mathcal{C}_{2}$ such that $$\mathfrak{o}(c_{2})\subseteq\cl_{j}(\mathfrak{o}(c_{2}))\subseteq\mathfrak{o}(a_{2})\subseteq \mathfrak{o}(c_{1})\cap\mathfrak{o}(x_{2}).$$  Since $\mathfrak{o}(x_{3})$ is $i$-dense and $\mathfrak{o}(c_{2})$ is a non-void $i$-open sublocale, it follows that $\mathfrak{o}(c_{2})\cap\mathfrak{o}(x_{3})\neq\mathsf{O}$. Applying that $(L,L_{1},L_{2})$ is $i$-prefit again implies that there is a nonzero $a_{3}\in L_{i}$ such that $$\mathsf{O}\neq\mathfrak{o}(a_{3})\subseteq\cl_{i}(\mathfrak{o}(a_{3}))\subseteq \mathfrak{o}(c_{2})\cap\mathfrak{o}(x_{3}).$$ Therefore, for the $i$-$\pi$-base $\mathcal{C}_{3}$, there is $\mathfrak{o}(c_{3})\in \mathcal{C}_{3}$ such that $$\mathfrak{o}(c_{3})\subseteq\cl_{j}(\mathfrak{o}(c_{3}))\subseteq\mathfrak{o}(a_{3})\subseteq \mathfrak{o}(c_{2})\cap\mathfrak{o}(x_{3}).$$ Continuing like this for $n=4,5,....$, we get that for each $i$-$\pi$-base $\mathcal{C}_{n}$, there is $\mathfrak{o}(c_{n})\in \mathcal{C}_{n}$ such that $$\mathfrak{o}(c_{n})\subseteq\cl_{j}(\mathfrak{o}(c_{n}))\subseteq \mathfrak{o}(c_{n-1})\cap\mathfrak{o}(x_{n}).$$ Since $(L,L_{1},L_{2})$ is $i$-pseudocomplete, $\bigcap_{n\in\mathbb{N}}\mathfrak{o}(c_{n})\neq\mathsf{O}$. Because $$\mathfrak{o}(c_{n})\subseteq \mathfrak{o}(y)\cap\mathfrak{o}(x_{n})$$ for each $n\in \mathbb{N}$, we have $$\mathsf{O}\neq \bigcap_{n\in\mathbb{N}}\mathfrak{o}(c_{n})\subseteq\mathfrak{o}(y)\cap\bigcap_{n\in\mathbb{N}}\mathfrak{o}(x_{n}),$$ making $\bigcap_{n\in\mathbb{N}}\mathfrak{o}(x_{n})$ an $i$-dense sublocale. Thus $(L,L_{1},L_{2})$ is $(i,j)$-Baire.
\end{proof}
Recall from \cite{BB} that the triple $(\mathfrak{J}L,(\mathfrak{J}L)_{1},(\mathfrak{J}L)_{2})$, where $\mathfrak{J}L$ is the locale of all ideals of $L$ and $(\mathfrak{J}L)_{i}$ ($i=1,2$) is the subframe of $\mathfrak{J}L$ consisting of all ideals $J\subseteq L$ generated by $J\cap L_{i}$, is a bilocale called the \textit{ideal bilocale}. 

Call a bilocale $(L,L_{1},L_{2})$ \textit{Noetherian} in case its total part $L$ is Noetherian, i.e., all of its elements are compact. In a Noetherian locale, all ideals are principal \cite{BFG}.

For use below, we recall from \cite[Proposition 6.9.]{N2} that in a bilocale $(L,L_{1},L_{2})$, if $x\in L_{i}$ is $j$-dense, then ${\downarrow}{x}\in \mathfrak{J}L_{i}$ is $\mathfrak{J}L_{j}$-dense. Furthermore, for a Noetherian bilocale $(L,L_{1},L_{2})$, $\bigvee\! J\in L_{i}$ is $j$-dense whenever $J\in \mathfrak{J}L_{i}$ is $\mathfrak{J}L_{j}$-dense.

\begin{proposition}
	A Noetherian bilocale $(L,L_{1},L_{2})$ is $(i,j)$-Baire iff $(\mathfrak{J}L,(\mathfrak{J}L)_{1},(\mathfrak{J}L)_{2})$ is $(i,j)$-Baire.
\end{proposition}
\begin{proof}
	$(\Longrightarrow)$: Let $\{\mathfrak{o}_{\mathfrak{J}L}(I_{n}):n\in\mathbb{N}\}$ be a collection of $\mathfrak{J}L_{i}$-dense $\mathfrak{J}L_{j}$-open sublocales. It follows that each $\bigvee I_{n}\in L_{j}$ is $i$-dense, making each $\mathfrak{o}\left(\bigvee I_{n}\right)$ an $i$-dense $j$-open sublocale of $L$. Pick a $\mathfrak{J}L_{i}$-open sublocale $\mathfrak{o}_{\mathfrak{J}L}(J)$ such that $$\mathfrak{o}_{\mathfrak{J}L}(J)\cap\left(\bigwedge_{n\in\mathbb{N}}\mathfrak{o}_{\mathfrak{J}L}(I_{n})\right)=\mathsf{O}.$$ 
	
	Claim: $\mathfrak{o}(\bigvee J)\cap \left(\bigwedge_{n\in\mathbb{N}}\mathfrak{o}(\bigvee I_{n})\right)=\mathsf{O}$.
	
	Proof: If not, then $$\mathsf{O}\neq\bigwedge_{n\in\mathbb{N}}\left(\mathfrak{o}\left(\bigvee J\right)\cap \mathfrak{o}\left(\bigvee I_{n}\right)\right)=\bigwedge_{n\in\mathbb{N}}\mathfrak{o}\left(\left(\bigvee J\right)\wedge \left(\bigvee I_{n}\right)\right)$$ so that $\left(\bigvee J\right)\wedge \left(\bigvee I_{n}\right)\neq0$ for each $n\in\mathbb{N}$. Therefore $$\mathsf{O}\neq{\downarrow}{\left(\bigvee J\right)}\cap {\downarrow}{\left(\bigvee I_{n}\right)}=J\cap I_{n}$$ where the latter equality follows since $(L,L_{1},L_{2})$ is Noetherian. We get that $\mathfrak{o}_{\mathfrak{J}L}(J)\cap \mathfrak{o}_{\mathfrak{J}L}(I_{n})\neq\mathsf{O}$ for each $n\in \mathbb{N}$, making 
	$$\mathfrak{o}_{\mathfrak{J}L}(J)\cap\left(\bigwedge_{n\in\mathbb{N}} \mathfrak{o}_{\mathfrak{J}L}(I_{n})\right)\neq\mathsf{O}$$ which is impossible.
	
	Now, since $(L,L_{1},L_{2})$ is $(i,j)$-Baire and each $\mathfrak{o}\left(\bigvee I_{n}\right)$ is $i$-dense $j$-open, it follows that $\bigwedge_{n\in\mathbb{N}}\mathfrak{o}(\bigvee I_{n})$ is $i$-dense so that $\mathfrak{o}(\bigvee J)=\mathsf{O}$. Therefore $\bigvee J=0$ which implies that $0={\downarrow}{\bigvee J}=J$. Thus $\mathfrak{o}_{\mathfrak{J}L}(J)=\mathsf{O}$. Hence $\bigwedge_{n\in\mathbb{N}}\mathfrak{o}_{\mathfrak{J}L}(I_{n})$ is $\mathfrak{J}L_{i}$-dense.
	
	$(\Longleftarrow)$: Choose a collection $\{\mathfrak{o}(x_{n}):n\in\mathbb{N}\}$ of $i$-dense $j$-open sublocales. Then $\{\mathfrak{o}_{\mathfrak{J}L}({\downarrow}{x_{n}}):n\in\mathbb{N}\}$ is a collection of $\mathfrak{J}L_{i}$-dense $\mathfrak{J}L_{j}$-open sublocales. By hypothesis, $\bigwedge_{n\in\mathbb{N}}\mathfrak{o}_{\mathfrak{J}L}({\downarrow}{x_{n}})$ is $\mathfrak{J}L_{i}$-dense. To show that $\bigwedge_{n\in\mathbb{N}}\mathfrak{o}(x_{n})$ is $i$-dense, let $\mathfrak{o}(y)$ be an $i$-open sublocale such that $\mathfrak{o}(y)\cap \left(\bigwedge_{n\in\mathbb{N}}\mathfrak{o}(x_{n})\right)=\mathsf{O}$. 
	
	Claim: $\mathfrak{o}_{\mathfrak{J}L}({\downarrow}{y})\cap \left(\bigwedge_{n\in\mathbb{N}}\mathfrak{o}_{\mathfrak{J}L}(\bigvee {\downarrow}{x_{n}})\right)=\mathsf{O}$.
	
	Proof: Otherwise, $\bigwedge_{n\in\mathbb{N}}\mathfrak{o}({\downarrow}{y}\cap{\downarrow}{x_{n}})\neq\mathsf{O}$ which implies that $$\mathsf{O}\neq {\downarrow}{y}\cap{\downarrow}{x_{n}}={\downarrow}{(y\wedge x_{n})}$$ for each $n\in\mathbb{N}$. Therefore $y\wedge x_{n}\neq 0$ for each $n\in\mathbb{N}$ so that $$\mathsf{O}\neq \bigwedge_{n\in\mathbb{N}}\mathfrak{o}(y\wedge x_{n})=\mathfrak{o}(y)\cap \left(\bigwedge_{n\in\mathbb{N}}\mathfrak{o}(x_{n})\right)$$ which is a contradiction.
	
	Thus $\mathfrak{o}_{\mathfrak{J}L}({\downarrow}{y})=\mathsf{O}$ implying that ${\downarrow}{y}=\mathsf{O}$. Therefore $\mathfrak{o}(y)=\mathsf{O}$ and hence $\bigwedge_{n\in\mathbb{N}}\mathfrak{o}(x_{n})$ is $i$-dense. 
\end{proof}
\begin{observation}
	In the proof of the reverse direction of the preceding result, we did not use the fact that $(L,L_{1},L_{2})$ is Noetherian.
\end{observation}

	%###############
%###############

\section{Concerning relative versions of $(i,j)$-Baire bilocales}

In this section, we consider $(i,j)$-Baireness of subbilocales.

We recall the following lemma from \cite{N1}.

\begin{lemma}\label{nudense}
	Let $(S,S_{1},S_{2})$ be a dense subbilocale of a bilocale $(L,L_{1},L_{2})$. An element $y$ of $L_{i}$ is $j$-dense iff $\nu_{S}(y)$ is $j_{S}$-dense.
\end{lemma}
\begin{corollary}\label{jsdenseopen}
	Let $(S,S_{1},S_{2})$ be a dense subbilocale of a bilocale $(L,L_{1},L_{2})$. An element $y$ of $L_{i}$ is $j$-dense iff $\mathfrak{o}_{S}(\nu_{S}(y))=S\cap\mathfrak{o}(y)$ is $j_{S}$-dense $i_{S}$-open.
\end{corollary}

We also have the following result.

\begin{lemma}\label{dense}
	Let $(L,L_{1},L_{2})$ be a bilocale with $(S,S_{1},S_{2})$ as its dense subbilocale. A sublocale $A$ of $S$ is $i_{S}$-dense iff it is $i$-dense.
\end{lemma}
\begin{proof}
	$(\Longrightarrow)$: Choose a non-void $i$-open sublocale $\mathfrak{o}(x)$ of $L$. Then $$\mathsf{O}\neq S\cap \mathfrak{o}(x)=\mathfrak{o}(\nu_{S}(x))$$ where $\nu_{S}(x)\in S_{i}$. This makes $\mathfrak{o}_{S}(\nu_{S}(x))$ an non-void $i_{S}$-open sublocale of $S$. Since $A$ is $i_{S}$-open, $$\mathsf{O}\neq A\cap\mathfrak{o}_{S}(\nu_{S}(x))=A\cap\mathfrak{o}(x).$$ Thus $A$ is $i$-dense.
	
	$(\Longleftarrow)$: Let $\mathfrak{o}_{S}(x)$ be a non-void $i_{S}$-open sublocale of $S$. Then $x=\nu_{S}(y)$ for some $y\in L_{i}$. It follows from Lemma \ref{nudense} that $y$ is $i$-dense. Therefore $$\mathsf{O}\neq A\cap\mathfrak{o}(y)=A\cap\mathfrak{o}_{S}(x).$$ Thus $A$ is $i_{S}$-dense.
\end{proof}
\begin{proposition}
	A bilocale $(L,L_{1},L_{2})$ is $(i,j)$-Baire only if it contains some dense $(i,j)$-Baire subbilocale.
\end{proposition}
\begin{proof}
	
	Let $(S,S_{1},S_{2})$ be a dense and $(i,j)$-Baire subbilocale of $(L,L_{1},L_{2})$ and pick a collection $\{\mathfrak{o}(x_{n}):n\in\mathbb{N}\}$ of $i$-dense $j$-open sublocales. Since the subbilocale $(S,S_{1},S_{2})$ is dense, it follows from Corollary \ref{jsdenseopen} that $\{S\cap\mathfrak{o}(x_{n}):n\in \mathbb{N}\}$ is a collection of $i_{S}$-dense $j_{S}$-open sublocales. By hypothesis, $\bigwedge_{n\in\mathbb{N}}(S\cap\mathfrak{o}(x_{n}))$ is $i_{S}$-dense, so that it is $i$-dense by Lemma \ref{dense}. Since $$\bigwedge_{n\in\mathbb{N}}(S\cap\mathfrak{o}(x_{n}))\subseteq\bigwedge_{n\in\mathbb{N}}\mathfrak{o}(x_{n}),$$ it follows that $\bigwedge_{n\in\mathbb{N}}\mathfrak{o}(x_{n})$ is $i$-dense.
\end{proof}
\begin{corollary}
		A bilocale $(L,L_{1},L_{2})$ is $(i,j)$-Baire only if $(\mathfrak{B}L,\nu_{\mathfrak{B}}[L_{1}],\nu_{\mathfrak{B}}[L_{2}])$ is $(i,j)$-Baire as a bilocale.
\end{corollary}
Call a bilocale $(L,L_{1},L_{2})$ \textit{$(i,j)$-submaximal} if every $i$-dense sublocale of $L$ is $j$-open
\begin{proposition}
	Let $(L,L_{1},L_{2})$ be an $(i,j)$-submaximal bilocale. Then $(L,L_{1},L_{2})$ is $(i,j)$-Baire iff $(\mathfrak{B}L,\nu_{\mathfrak{B}}[L_{1}],\nu_{\mathfrak{B}}[L_{2}])$ is $(i,j)$-Baire as a bilocale.
\end{proposition}
\begin{proof}
	We only prove the forward implication:
	
	Let $\{\mathfrak{o}_{\mathfrak{B}L}(x_{n}):n\in\mathbb{N}\}$ be a collection of $i_{\mathfrak{B}L}$-dense $j_{\mathfrak{B}L}$-open sublocales. It follows that $\{\mathfrak{o}(x_{n}):n\in\mathbb{N}\}$ is a collection of $i$-dense $j$-open sublocales. Since $(L,L_{1},L_{2})$ is $(i,j)$-Baire, $\bigwedge_{n\in\mathbb{N}}\mathfrak{o}(x_{n})$ is $i$-dense. We must have that  $\bigwedge_{n\in\mathbb{N}}\mathfrak{o}_{\mathfrak{B}L}(x_{n})$ is $i_{\mathfrak{B}L}$-dense, otherwise there is a non-void $\nu_{\mathfrak{B}L}[L_{i}]$-open sublocale $\mathfrak{o}_{\mathfrak{B}L}(y)$ such that $$\mathfrak{o}_{\mathfrak{B}L}(y)\cap \left(\bigwedge_{n\in\mathbb{N}}\mathfrak{o}_{\mathfrak{B}L}(x_{n})\right)=\mathsf{O}.$$ Therefore $$\mathfrak{o}_{\mathfrak{B}L}(y)\cap \left(\bigwedge_{n\in\mathbb{N}}\mathfrak{o}(x_{n})\right)=\mathsf{O}.$$ Since every dense sublocale is $i$-dense and $(L,L_{1},L_{2})$ is $(i,j)$-submaximal, we have that $\mathfrak{B}L$ is $j$-open so that $\mathfrak{o}_{\mathfrak{B}L}(y)=\mathfrak{B}L\cap\mathfrak{o}(y)$ is a $j$-open sublocale.  Therefore $\mathfrak{o}_{\mathfrak{B}L}(y)=0$ which is impossible.
\end{proof}

\begin{proposition}
	Every $i$-open subbilocale of an $(i,j)$-Baire bilocale is $(i,j)$-Baire.
\end{proposition}
\begin{proof}
	Let $(S,S_{1},S_{2})$ be an $i$-open subbilocale of an $(i,j)$-Baire bilocale $(L,L_{1},L_{2})$. Choose a collection $\{\mathfrak{o}_{S}(x_{n}):n\in\mathbb{N}\}$ of $i_{S}$-dense $j_{S}$-open sublocales. We show that $\bigwedge_{n\in\mathbb{N}}\mathfrak{o}_{S}(x_{n})$ is $i_{S}$-dense. Pick an $i_{S}$-open sublocale $\mathfrak{o}_{S}(y)$ such that $$\left(\bigwedge_{n\in\mathbb{N}}\mathfrak{o}_{S}(x_{n})\right)\cap \mathfrak{o}_{S}(y)=\mathsf{O}.$$ Since $\mathfrak{o}_{S}(y)\subseteq\overline{S}$, $\mathfrak{o}_{S}(y)\cap (L\smallsetminus\overline{S})=\mathsf{O}$. Therefore \begin{align*}
		\mathfrak{o}_{S}(y)\cap\left(\left(\bigwedge_{n\in\mathbb{N}}\mathfrak{o}_{S}(x_{n})\right)\vee (L\smallsetminus\overline{S})\right)&=\left(\left(\bigwedge_{n\in\mathbb{N}}\mathfrak{o}_{S}(x_{n})\right)\cap \mathfrak{o}_{S}(y)\right)\vee \left(\mathfrak{o}_{S}(y)\cap (L\smallsetminus\overline{S})\right)\\
		&=\left(\bigwedge_{n\in\mathbb{N}}\mathfrak{o}_{S}(x_{n})\right)\cap \mathfrak{o}_{S}(y)\\
		&=\mathsf{O}.
	\end{align*}
Because $\mathfrak{o}_{S}(x_{n})\vee (L\smallsetminus\overline{ S})$ is $i$-dense, it follows that $$\bigwedge_{n\in\mathbb{N}}\left(\mathfrak{o}_{S}(x_{n})\vee(L\smallsetminus\overline{ S})\right)=(L\smallsetminus\overline{ S})\vee \bigwedge_{n\in\mathbb{N}}\mathfrak{o}_{S}(x_{n})$$ is $i$-dense. Therefore $\mathfrak{o}_{S}(y)=\mathsf{O}$. Thus $\bigwedge_{n\in\mathbb{N}}\mathfrak{o}_{S}(x_{n})$ is $i_{S}$-dense.
\end{proof}

\begin{definition}
	Let $(L,L_{1},L_{2})$ be a bilocale. A subbilocale $(S,S_{1},S_{2})$ of $(L,L_{1},L_{2})$ is \textit{relatively $(i,j)$-Baire} if for every collection $\{\mathfrak{o}(x_{n}):n\in\mathbb{N}\}$ of $i$-dense $j$-open sublocales, $S\cap \left(\bigwedge_{n\in\mathbb{N}}\mathfrak{o}(x_{n})\right)$ is $i_{S}$-dense.
\end{definition}

\begin{proposition}
	In a class of dense subbilocales, $(i,j)$-Baire coincides with relatively $(i,j)$-Baire.
\end{proposition}
\begin{proof}
	Let $(S,S_{1},S_{2})$ be an $(i,j)$-Baire subbilocale of a bilocale $(L,L_{1},L_{2})$ and choose a collection $\{\mathfrak{o}(x_{n}):n\in\mathbb{N}\}$ of $i$-dense $j$-open sublocales of $L$. If $$\mathfrak{o}_{S}(y)\cap S\cap\left(\bigwedge_{n\in\mathbb{N}}\mathfrak{o}(x_{n})\right)=\mathsf{O},$$ then $$\mathsf{O}=\mathfrak{o}_{S}(y)\cap \left(\bigwedge_{n\in\mathbb{N}}(S\cap\mathfrak{o}(x_{n}))\right)=\mathfrak{o}_{S}(y)\cap\left(\bigwedge_{n\in\mathbb{N}}\mathfrak{o}_{S}(\nu_{S}(x_{n}))\right)$$ where each $\mathfrak{o}_{S}(\nu_{S}(x_{n}))$ is $i_{S}$-dense and $j_{S}$-open. Since $(S,S_{1},S_{2})$ is $(i,j)$-Baire as a bilocale, $\bigwedge_{n\in\mathbb{N}}\mathfrak{o}_{S}(\nu_{S}(x_{n}))$ is $i_{S}$-dense so that $\mathfrak{o}_{S}(y)=\mathsf{O}$. Thus $S\cap\left(\bigwedge_{n\in\mathbb{N}}\mathfrak{o}(x_{n})\right)$ is $i_{S}$-dense. 
	
	On the other hand, let $(S,S_{1},S_{2})$ be a relatively $(i,j)$-Baire subbilocale and pick a collection  $\{\mathfrak{o}_{S}(x_{n}):n\in\mathbb{N}\}$ of $i_{S}$-dense $j_{S}$-open sublocales of $S$. For each $x_{n}$, there is $a_{n}\in L_{j}$ such that $x_{n}=\nu_{S}(a_{n})$. Now, members of the collection $\{\mathfrak{o}_{S}(a_{n}):n\in\mathbb{N}\}$ are $i$-dense $j$-open in $(L,L_{1},L_{2})$. Since $(S,S_{1},S_{2})$ is relatively $(i,j)$-Baire, $$S\cap\left(\bigwedge_{n\in\mathbb{N}}\mathfrak{o}(a_{n})\right)=\bigwedge_{n\in\mathbb{N}}\mathfrak{o}_{S}(x_{n})$$ is $i_{S}$-dense. Thus $(S,S_{1},S_{2})$ is $(i,j)$-Baire.
\end{proof}

We close this section with a characterization of relatively $(i,j)$-Baire subbilocales.
\begin{proposition}
	Let $(S,S_{1},S_{2})$ be a dense and complemented subbilocale of a bilocale $(L,L_{1},L_{2})$ whose $j$-G$_{\delta}$-sublocales are complemented. The following statements are equivalent:
	\begin{enumerate}
		\item $(S,S_{1},S_{2})$ is relatively $(i,j)$-Baire.
		\item For every non-void $i$-open sublocale $U$ of $L$, $S\cap U$ is of $(j,i)$-second category in $(S,S_{1},S_{2})$.
		\item For every sublocale $U$ of $(j,i)$-first category in $(L,L_{1},L_{2})$, $\inte_{i_{S}}(S\cap U)=\mathsf{O}$. 
		\item If $V$ is a sublocale of $(j,i)$-first category in $(L,L_{1},L_{2})$, then  $S\cap (L\smallsetminus V)$ is $i_{S}$-dense. 
	\end{enumerate} 
\end{proposition}
\begin{proof}
		$(1)\Longrightarrow(2)$: Let $\mathfrak{o}(x)$ be non-void $i$-open and assume that $$S\cap\mathfrak{o}(x)\subseteq\bigvee_{n\in\mathbb{N}}^{S}\mathfrak{c}_{S}(x_{n})$$ for some collection $\{\mathfrak{c}_{S}(x_{n}):n\in\mathbb{N}\}$ of $(j_{S},i_{S})$-nowhere dense sublocales of $S$. Then $$S\cap\mathfrak{o}(x)\subseteq\bigvee_{n\in\mathbb{N}}\mathfrak{c}(x_{n})$$ where each $\mathfrak{c}(x_{n})$ is $(j,i)$-nowhere dense because $(S,S_{1},S_{2})$ is dense. It is clear that the collection $\{\mathfrak{o}(x_{n}):n\in\mathbb{N}\}$ consists of $i$-dense $j$-open sublocales. It follows from $(1)$ that $S\cap\left(\bigwedge_{n\in\mathbb{N}}\mathfrak{o}(x_{n})\right)$ is $i_{S}$-dense. Since $S\cap\mathfrak{o}(x)\neq\mathsf{O}$ because of density of $(S,S_{1},S_{2})$, we have that $$S\cap\mathfrak{o}(x)\cap S\cap\left(\bigwedge_{n\in\mathbb{N}}\mathfrak{o}(x_{n})\right)\neq\mathsf{O}.$$Therefore $$\left(\bigvee_{k\in\mathbb{N}}\mathfrak{c}(x_{k})\right)\cap \left(\bigwedge_{n\in\mathbb{N}}\mathfrak{o}(x_{n})\right)\neq\mathsf{O}.$$ Since $\bigwedge_{n\in\mathbb{N}}\mathfrak{o}(x_{n})$ is a $j$-G$_{\delta}$-sublocale of $L$, it follows that it is complemented. Therefore $$\mathsf{O}\neq\bigvee_{k\in\mathbb{N}}\left(\mathfrak{c}(x_{k})\cap \left(\bigwedge_{n\in\mathbb{N}}\mathfrak{o}(x_{n})\right)\right)\subseteq\bigvee_{k\in\mathbb{N}}\left(\mathfrak{c}(x_{k})\cap \mathfrak{o}(x_{k})\right)=\bigvee_{n\in\mathbb{N}}(\mathsf{O})=\mathsf{O}$$which is impossible. Thus $S\cap \mathfrak{o}(x)$ is $(j,i)$-second category.
		
		$(2)\Longrightarrow(1)$: Let $\{\mathfrak{o}(x_{n}):n\in\mathbb{N}\}$ be a collection of $i$-dense $j$-open sublocales and assume that there is non-void $i_{S}$-open sublocale $\mathfrak{o}_{S}(y)$ of $S$ such that $$\mathfrak{o}_{S}(y)\cap S\cap\left(\bigwedge_{n\in\mathbb{N}}\mathfrak{o}(x_{n})\right)=\mathsf{O}.$$ Then $\mathfrak{o}(y)$ is non-void $i$-open and $$\mathfrak{o}_{S}(y)\cap\left(\bigwedge_{n\in\mathbb{N}}\mathfrak{o}(x_{n})\right)=\mathsf{O}$$ which implies $$\mathfrak{o}_{S}(y)\subseteq S\cap\left(\bigvee_{n\in\mathbb{N}}\mathfrak{o}(x_{n})\right)=\bigvee_{n\in\mathbb{N}}\mathfrak{c}_{S}(\nu_{S}(x_{n}))$$ where the latter equality holds since $S$ is complemented. Since each $\mathfrak{c}_{S}(\nu_{S}(x_{n}))$ is $(j_{S},i_{S})$-nowhere dense, $S\cap\mathfrak{o}(y)=\mathfrak{o}_{S}(y)$ is of $(j,i)$-first category in $(S,S_{1},S_{2})$ which is a contradiction.
		
		$(2)\Longrightarrow(3)$: Let $U$ be a sublocale of $L$ which is of $(j,i)$-first category in $(L,L_{1},L_{2})$ and assume that $\inte_{i_{S}}(S\cap U)\neq\mathsf{O}$. Then $$\inte_{i_{S}}(S\cap U)=\mathfrak{o}(x)\cap S$$ for some $x\in L_{i}$. Such $\mathfrak{o}(x)$ is a non-void $i$-open sublocale of $L$, so $$\inte_{i_{S}}(S\cap U)=S\cap\mathfrak{o}(x)$$ must be of $(j,i)$-second category in $(S,S_{1},S_{2})$ by (2). But $U\subseteq\bigvee_{n\in\mathbb{N}}\mathfrak{c}(x_{n})$ for some collection $\{\mathfrak{c}(x_{n}):n\in\mathbb{N}\}$ of $(j,i)$-nowhere dense sublocales of $L$, so $$\inte_{i_{S}}(S\cap U)=\mathfrak{o}(x)\cap S\subseteq U\cap S\subseteq S\cap \bigvee_{n\in\mathbb{N}}\mathfrak{c}(x_{n})=\bigvee_{n\in\mathbb{N}}\mathfrak{c}_{S}(\nu_{S}(x_{n}))$$ where each $\mathfrak{c}_{S}(\nu_{S}(x_{n}))$ is $(j_{S},i_{S})$-nowhere dense in $(S,S_{1},S_{2})$. This makes $\mathfrak{o}(x)\cap S$ a sublocale of $(j,i)$-first category in $(S,S_{1},S_{2})$ which is impossible.
		
		$(3)\Longrightarrow(4)$: Let $V$ be a sublocale of $L$ which is of $(j,i)$-first category in $(L,L_{1},L_{2})$ and choose an $i_{S}$-open sublocale $\mathfrak{o}_{S}(x)$ such that $$\mathfrak{o}_{S}(x)\cap S\cap(L\smallsetminus V)=\mathsf{O}.$$ Then $$\mathfrak{o}_{S}(x)\subseteq S\cap V.$$ By (3), $\inte_{ji{S}}(S\cap V)=\mathsf{O}$, making $\mathfrak{o}_{S}(x)=\mathsf{O}$. 
		
		$(4)\Longrightarrow(2)$: Let $\mathfrak{o}(x)$ be a non-void $i$-open sublocale of $L$ and assume that $S\cap\mathfrak{o}(x)$ is of $(j,i)$-first category. By $(4)$, $$S\cap(L\smallsetminus \mathfrak{o}(x))=S\cap\mathfrak{c}(x)=\mathfrak{c}_{S}(\nu_{S}(x))$$ is $i_{S}$-dense which implies that $\nu_{S}(x)=0$. Therefore $\mathfrak{o}(x)=\mathsf{O}$ which is a contradiction.
\end{proof}

\section{Baireness of topobilocales}

The aim of this section is to introduce and characterize Baireness in the category of topobilocales. 

A \textit{topobilocale} \cite{ZEFB} is a triple $(L,\tau_{1},\tau_{2})$ where $L$ is a locale, $L_{1}$ and $L_{2}$ are subframes of $L$ all of whose elements are complemented in $L$. Each member of $\tau_{i}$ ($i=1,2$) is called $\tau_{i}$-open.

For each $a\in L$, the \textit{$\tau_{i}$-closure} ($i=1,2$) of $a$ in $L$ is defined by $$\cl_{(L,\tau_{i})}(a)=\bigwedge\{b\in \tau_{i}':a\leq b\}$$ and the \textit{$\tau_{i}$-interior} of $a$ is defined by $$\inte_{(L,\tau_{i})}(a)=\bigvee\{b\in \tau_{i}:b\leq a\}.$$

We have the following result. See \cite{ZE} for the proofs of some of the statements. For the rest of the statements, the proofs of the other statements resemble that of  \cite[Proposition 5.1.3.]{N1}. 
\begin{proposition}
	Let $(L,\tau_{i},\tau_{j})$ be a topobilocale. Then
	\begin{enumerate}
		\item $\cl_{(L,\tau_{i})}(0)=\inte_{(L,\tau_{i})}(0)=0$.
		\item $\cl_{(L,\tau_{i})}(1)=\inte_{(L,\tau_{i})}(1)=1$.
			\item $a\leq\cl_{(L,\tau_{i})}(a)$.
			\item If $a\leq b$, then $\cl_{(L,\tau_{i})}(a)\leq \cl_{(L,\tau_{i})}(b)$.
			\item $\inte_{(L,\tau_{i})}(a)\leq a$.
			\item If $a\leq b$, then $\inte_{(L,\tau_{i})}(a)\leq\inte_{(L,\tau_{i})}(b)$.
			\item For each $a\in L$, $(\cl_{(L,\tau_{i})}(a))'=\inte_{(L,\tau_{i})}(a')$.
			\item For each $a\in L$, $(\inte_{(L,\tau_{i})}(a))'=\cl_{(L,\tau_{i})}(a')$.
	\end{enumerate}
\end{proposition}

Call an element $a\in L$ \textit{$\tau_{i}$-dense} if $\cl_{(L,\tau_{i})}(a)=1$.

\begin{definition}
	Call a topobilocale $(L,\tau_{1},\tau_{2})$ \textit{$(\tau_{i},\tau_{j})$-Baire} if any sequence $(x_{n})_{n\in\mathbb{N}}$ of $\tau_{i}$-dense elements of $\tau_{j}$ satisfies the condition $\bigwedge_{n\in\mathbb{N}}x_{n}$ is $\tau_{i}$-dense. 
\end{definition}

Call an element $a\in L$ \textit{$(\tau_{i},\tau_{j})$-nowhere dense} if $\inte_{(L,\tau_{j})}(\cl_{(L,\tau_{i})}(a))=0$.

An element $a\in L$ is of \textit{$(\tau_{i},\tau_{j})$-first category} if $a\leq \bigvee_{n\in\mathbb{N}}x_{n}$ for some collection $\{x_{n}:n\in\mathbb{N}\}$ of $(\tau_{i},\tau_{j})$-nowhere dense elements of $L$. Otherwise it is of \textit{$(\tau_{i},\tau_{j})$-second category}.

For use below, we give the following result with a proof similar to that of \cite[Proposition 2.1.4.]{N2}

\begin{proposition}
	Let $(L,\tau_{1},\tau_{2})$ be a topobilocale. Then $a\in L$ is $(\tau_{i},\tau_{j})$-nowhere dense iff $(\cl_{(L,\tau_{i})}(a))'$ is $\tau_{j}$-dense.
\end{proposition}

 The proof of the following result is similar to that of Proposition \ref{mainresult}. We only prove $(1)\Longrightarrow(2)$.
\begin{proposition}
	Let $(L,\tau_{1},\tau_{2})$ be a topobilocale. The following statements are equivalent.
	\begin{enumerate}
		\item $(L,\tau_{1},\tau_{2})$ is $(\tau_{i},\tau_{j})$-Baire.
		\item Each nonzero $\tau_{i}$ element is of $(\tau_{j},\tau_{i})$-second category.
		\item Every element of $(\tau_{j},\tau_{i})$-first category has a zero $\tau_{i}$-interior.
		\item The complement an element of $(\tau_{j},\tau_{i})$-first category is $\tau_{i}$-dense. 
	\end{enumerate}
\end{proposition}
\begin{proof}
	$(1)\Longrightarrow(2)$: Assume that there is a nonzero element $a\in \tau_{i}$ which is of $(\tau_{j},\tau_{i})$-first category. Then $a\leq \bigvee_{n\in\mathbb{N}}x_{n}$ for some collection $\{x_{n}:n\in\mathbb{N}\}$ of $(\tau_{j},\tau_{i})$-nowhere dense elements. It is clear members of the collection $\{(\cl_{(L,\tau_{j})}(x_{n}))':n\in\mathbb{N}\}$ are $\tau_{i}$-dense. By (1), $\bigwedge_{n\in\mathbb{N}}(\cl_{(L,\tau_{j})}(x_{n}))'$ is $\tau_{i}$-dense. It follows that $$a\wedge\left(\bigwedge_{n\in\mathbb{N}}(\cl_{(L,\tau_{j})}(x_{n}))'\right)\neq 0.$$ Therefore \begin{align*}
		0&\neq\left(\bigvee_{k\in\mathbb{N}}x_{k}\right)\wedge \left(\bigwedge_{n\in\mathbb{N}}(\cl_{(L,\tau_{j})}(x_{n}))'\right)\\
		&=\bigvee_{k\in\mathbb{N}}\left(x_{k}\wedge \left(\bigwedge_{n\in\mathbb{N}}(\cl_{(L,\tau_{j})}(x_{n}))'\right)\right)\text{ since }L \text{ is a locale}\\
		&\leq \bigvee_{k\in\mathbb{N}}\left(x_{k}\wedge (\cl_{(L,\tau_{j})}(x_{k}))'\right)\\
		&\leq \bigvee_{k\in\mathbb{N}}\left(\cl_{(L,\tau_{i})}(x_{k})\wedge \cl_{(L,\tau_{j})}(x_{k})\right)\\
		&=0
	\end{align*}which is a contradiction.
\end{proof}

	%############
	%############

\end{document}